\documentclass[12pt]{article}
\usepackage{amsmath,amssymb,amsthm,mathrsfs,graphicx,floatflt,color,epsfig,verbatim,multirow}
\newtheorem{thm}{Theorem}[section]
\newtheorem{prop}[thm]{Proposition}
\newtheorem{lemma}[thm]{Lemma}
\newtheorem{cor}[thm]{Corollary}
\theoremstyle{definition}
\newtheorem{ex}{Example}

\begin{document}

\title{Bijections Between Multiparking Functions, Dirichlet Configurations, and Descending $R$-Traversals}

\author{Dimitrije N. Kosti\'{c}\thanks{Partially supported by NSF VIGRE grant \# 9977354.}
\vspace{.3cm} \\
Department of Mathematics \\
Texas A\&M University,  College Station, TX 77843\\
dkostic@math.tamu.edu
\vspace{.4cm}
}
\date{}
\maketitle



\begin{abstract}

There are several combinatorial objects that are known to be in bijection to the spanning trees of a graph $G$.  These objects include $G$-parking 
functions, critical configurations of $G$, and descending traversals of $G$.  In this paper, we extend the bijections to generalizations of all three 
objects.

\end{abstract}


\begin{section}{Introduction}

A \emph{parking function} (of length $n$) is a sequence $(a_1, a_2, \ldots, a_n)$ of nonnegative integers for which there exists a permutation $\pi 
\in S_n$ so that $a_{\pi(i)} < i$ for every $1 \leq i \leq n$.  This concept was introduced (by way of an analogy to parking on a one-way street) in 
a 1966 paper by Konheim and Weiss \cite{KoWe} on linear probing.  In combinatorics, parking functions are most famous for being in bijection to the 
set of labelled trees on $n$ vertices and several combinatorial proofs of this fact are known; see \cite{FoRi} for an example.  More recently, 
renewed interest in parking functions has spurred various generalizations.  In \cite{Stan}, Stanley showed that $k$-parking functions are in 
bijection to the chambers of the extended Shi arrangements and in \cite{Yan} Yan showed that $k$-parking functions are in bijection to sequences of 
rooted $b$-forests.  In \cite{StPi}, Stanley and Pittman showed that the number of $\overrightarrow{x}$-parking functions is determined by the volume 
of a certain polytope.

Let $G$ be a loopless graph with vertex set $V(G) = \{ 1, 2, \ldots, n \}$. (Unless otherwise stated, we do not assume connectedness.) In 2004 
Postnikov and Shapiro \cite{PoSh} proposed a new generalization of the notion of a parking function, equivalent to the following.  Let $\mathbb{N}$ 
be the set of nonnegative integers and let $G$ be connected.  If $i \in U \subseteq V(G)$, the \emph{out-degree} of $i$ in $U$ is $\mathbb{O}_U(i) := 
\# \{ j \in V(G) - U \; | \; i$ is adjacent to $j \}$.  The \emph{in-degree} of $i$ in $U$ is $\mathbb{I}_U(i) := \# \{ j \in U \;
| \; i$ is adjacent to $j \}$.  A {\em $G$-parking function} is a function $f : V(G) \rightarrow \mathbb{N} \cup \{ \infty \}$ such that for any $U 
\subseteq V(G) - \{ 1 \}$ there exists a vertex $i \in U$ such that $0 \leq f(i) < \mathbb{O}_U(i)$.  In \cite{ChPy}, a family of bijections between 
the set of $G$-parking functions and spanning trees of $G$ was constructed, each bijection being determined by a {\em proper set of tree orders}.

Kosti\'{c} and Yan \cite{KoYa} generalized this work further.  A {\em $G$-multiparking function} is a function $f : V(G) \rightarrow \mathbb{N} \cup 
\{ \infty \}$ such that for any $U \subseteq V(G)$ there exists $i \in U$ with either {\bf (A)} $f(i) = \infty$, or {\bf (B)} $0 \leq f(i) < 
\mathbb{O}_U(i)$.  Throughout this paper, we will refer to those vertices $i$ with $f(i) = \infty$ as \emph{roots} and those with $0 \leq f(i) < 
\mathbb{O}_U(i)$ as being \emph{well-behaved in $U$}.  Let $\mathcal{MP} = \mathcal{MP}_{R,G}$ denote the set of $G$-multiparking functions with root 
set $R$.

There is a subtle but important difference between this definition of a $G$-multiparking function and the one that appears in \cite{KoYa}.  In that 
paper, the minimal vertex in each component of $G$ is required to be a root; here there is no such restriction.  Note, however, that $R$ cannot be 
empty; $V(G)$ cannot have a well-behaved vertex, so it must have a root.

This paper will construct bijections between $G$-multiparking functions and two other objects, Dirichlet configurations and descending 
$R$-traversals.  In section $2$, we provide an algorithm to verify whether a function is a $G$-multiparking function.  In section $3$ we introduce 
Dirichlet configurations and some basic facts about them.  In section $4$ we establish a bijection between $G$-multiparking functions and Dirichlet 
configurations on $G$.  In section $5$ we introduce descending $R$-traversals, define a certain partition of the set of these objects, and prove that 
they are in bijection to $G$-multiparking functions.

\end{section}


\begin{section}{A Burning Algorithm for Multiparking Functions}

We begin with a simple result on $G$-multiparking functions.  It is proven in \cite{KoYa}, but we present the proof here for completeness.

\begin{lemma}~\label{burning} Let $f : V(G) \rightarrow \mathbb{N} \cup \{ \infty \}$ be a vertex function.  Then, $f$ is a $G$-multiparking function 
if and only if there exists a permutation $\pi \in S_n$ such that the vertex $\pi(i)$ is either a root of $f$ or well-behaved in the set $U_i := V(G) 
- \{ \pi(1), \pi(2), \ldots, \pi(i-1) \}$.

\end{lemma}

\begin{proof} 

Suppose $\pi \in S_n$ satisfies the condition and let $U \subseteq V(G)$.  Let $l$ be the maximum index such that $U \subseteq U_l$.  By maximality, 
$\pi(l) \in U$.  Thus for $i=\pi(l)$, either $i$ is a root or $0 \leq f(i) < \mathbb{O}_{U_l}(i) \leq \mathbb{O}_U(i)$.  The other implication 
follows from the definition of a $G$-multiparking function. \end{proof}

This lemma suggests a simple algorithm for determining whether a function is a $G$-multiparking function.  Let $U_0=V(G)$.  If there is a vertex $i_1 
\in U_0$ that is a root (clearly, nothing in $U_0$ can be well-behaved), then let $U_1 := U_0 - \{ i_1 \}$.  If there is a vertex $i_2 \in U_1$ which 
is a root or well-behaved, then let $U_2 := U_1 - \{ i_2 \}$.  Lemma \ref{burning} implies that, when this process is continued, $U_n = \emptyset$ if 
and only if $f$ is a $G$-multiparking function.  This algorithm is a generalization of the \emph{burning algorithm}, which was originally developed 
in \cite{Dhar} to study critical configurations.

\begin{ex} In several examples in this paper, the following graph $\Gamma$ will be considered.

\begin{figure}[h]
\begin{picture}(0,0)
\put(200,-80){\includegraphics[width=4cm]{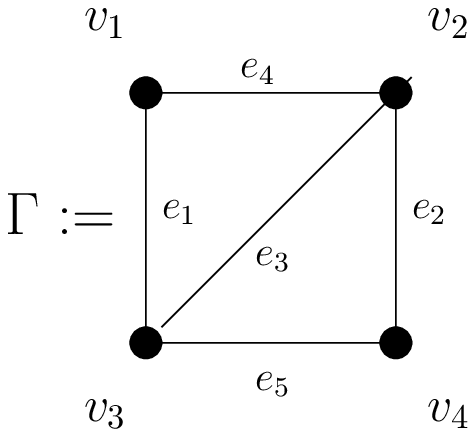}}
\end{picture}
\vskip 0.9 in
\end{figure}

The following example illustrates the burning algorithm described above, proving that the illustrated function is a $\Gamma$-(multi)parking function.  
In the leftmost picture the circled vertices are in $U_1$.  In the next leftmost picture the circled vertices are in $U_2$, and so forth.

\begin{figure}[h]
\begin{picture}(0,0)
\put(5,-95){\includegraphics[width=16.5cm]{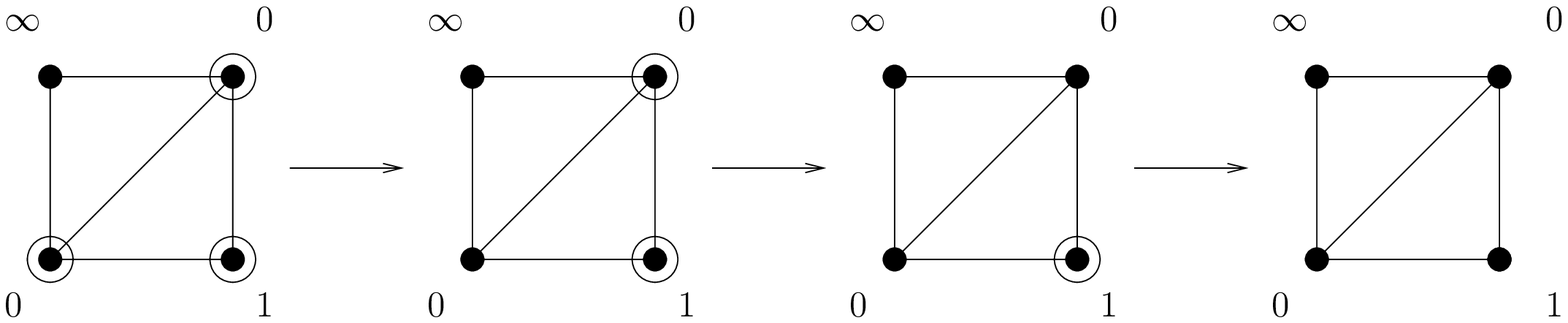}}
\end{picture}
\vskip 0.9 in
\end{figure}

\end{ex}

\end{section}

\begin{section}{Dirichlet Configurations}

Now we introduce a generalization of a structure that appears in the literature in a variety of contexts and, with minor variations, is known as a 
critical configuration, a sandpile model, and a chip-firing game.  Let $R$ be a set of vertices containing at least one vertex from each component 
of $G$.  A \emph{configuration} $\mu$ on $G$ (with \emph{root set} $R$) is an integer-valued function on the vertex set for which $\mu(i) = -\infty$ 
if $i \in R$ and $0 \leq \mu(i) < \infty$ otherwise.  A vertex $i$ is said to be \emph{ready} (in $\mu$) if $\mu(i) \geq deg(i)$.  $\mu$ is 
\emph{stable} if $0 \leq \mu(i) < deg(i)$ for every $i \notin R$.  An \emph{avalanche} is a finite sequence $\alpha = (\mu_1, \mu_2, \ldots \mu_t)$ 
of configurations on $G$, where for each $1 \leq s < t$ there exists a vertex $i_s \in V(G)$ which is ready in $\mu_s$ and
\[ \mu_{s+1}(i) = 
\left\{ \begin{array}{ll} 
\mu_s(i)-deg(i) & \mbox{if $i=i_s$} \\ 
\mu_s(i)+e(i,i_s) & \mbox{if $i \neq i_s$} 
\end{array} \right. \] 
\noindent where $e(i,i_s)$ is the number of edges between $i$ and $i_s$.  If we think of $\mu_s$ as keeping track of how many ``chips'' are stored 
at each vertex on the graph, then we transform $\mu_s$ into $\mu_{s+1}$ by sending a chip down each edge adjacent to $i_s$.  This process is often 
called \emph{firing} a vertex (hence the ``chip-firing'' terminology), and so one usually thinks of an avalanche as a sequence of vertex firings.  
Note that only vertices that are ready can be fired, and that the same vertex may be fired several times in succession if it has a large enough 
number of chips.  We say that $\alpha$ \emph{begins} at $\mu_1$, \emph{ends} at $\mu_t$, and \emph{connects} these two configurations.  We use the 
convention that if, in any avalanche, $\mu_1$ is stable then every vertex in $R$ is fired in some arbitrary but fixed order and that this is the 
only situation in which roots are fired.  Note that if a chip is sent to a root it disappears from the system; it follows from the connectedness of 
each component of $G$ that, given any configuration, there is an avalanche leading to a stable configuration.  $\mu$ is \emph{recurrent} if there is 
an avalanche that begins and ends at $\mu$.  $\mu$ is \emph{Dirichlet} if it is both stable and recurrent.  Let $\mathcal{DC} = \mathcal{DC}_{R,G}$ 
denote the set of Dirichlet configurations on $G$ with root set $R$.

Dirichlet configurations are usually called \emph{critical configurations} when $G$ is connected, and this case has been studied extensively (see, 
for example, \cite{Bigg1}).  Aspects of Dirichlet configurations were first examined in \cite{ChEl}, such as bounds on the number of vertex firings 
necessary to reach a stable configuration.

\begin{ex} The following example illustrates a critical configuration for $\Gamma$.  Every vertex is labelled ``$v_i / n$", where $v_i$ is the vertex 
label and $n$ is the number of chips at that vertex at that configuration.  The vertex about to be fired in each configuration is circled.

\begin{figure}[h]
\begin{picture}(0,0)
\put(5,-70){\includegraphics[width=16.5cm]{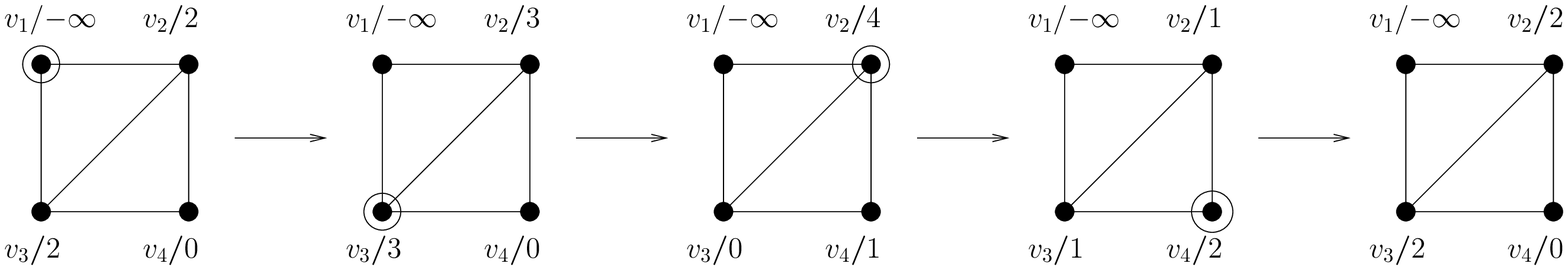}}
\end{picture}
\vskip 0.8 in
\end{figure}

\end{ex}

Note that in the second configuration in the above avalanche, we could have fired either $v_2$ or $v_3$.  If we had fired $v_2$ 
instead of $v_3$, we would still have ended the avalanche on the configuration we started with.  

It is not difficult to see that Dirichlet configurations exist on every graph (for instance, the configuration with $R=V(G)$) and that for every 
configuration there is an avalanche ending on a stable configuration.  This is essentially because every component contains a root and therefore the 
total number of chips on the graph is nonincreasing after the root firings (if the first configuration is stable) in an avalanche.  See Lemma 1 of 
\cite{ChEl} for a detailed proof.

The following is a characterization of recurrent configurations.  Let $\chi$ be the configuration

\[
\chi(i) =
\left\{ \begin{array}{ll}
0 & \mbox{if $i \in R$} \\
\sum_{r \in R} e(i,r) & \mbox{if $i \notin R$}
\end{array} \right.
\]

\begin{prop}~\label{recurrentcharacterization} The configuration $\mu$ is Dirichlet $\Longleftrightarrow$ it is stable and there is an avalanche 
connecting $(\mu + \chi)$ to $\mu$, where $(\mu + \chi)(v) = \mu(v) + \chi(v)$.

\end{prop}

\begin{proof} 

($\Longleftarrow$) The trivial avalanche (consisting of firing all the roots only) connects $\mu$ to $(\mu + \chi)$, and thus concatenating this 
avalanche with the avalanche connecting $(\mu + \chi)$ to $\mu$ shows that $\mu$ is recurrent.  Since it is stable, $\mu$ is Dirichlet.

($\Longrightarrow$) Given a Dirichlet $\mu$, it is stable and recurrent.  Thus there is an avalanche \\ $(\mu, \omega_1, \omega_2, \ldots, \omega_l, 
\mu)$.  But since $\mu$ is stable, the roots are the only vertices that can be fired first.  Thus, $\omega_k = (\mu + \chi)$ where $k$ is the number 
of roots.  Thus, $(\omega_k, \omega_{k+1}, \ldots, \omega_l, \mu)$ is the necessary avalanche.  \end{proof}

Cori and Rossin \cite{CoRo} have a similar proof for the case when the graph is connected.  The set of critical configurations of $G$ is closely 
related to the set of $G$-parking functions; the most famous connection is that both sets are in bijection to the spanning trees of $G$.  Here, 
however, we provide a bijection between $G$-multiparking functions and Dirichlet configurations on $G$ that does not go through the set of spanning 
trees.  To simplify the presentation, we will assume $G$ has no multiple edges.

\end{section}


\begin{section}{A Bijection Between Dirichlet Configurations and $G$-Multiparking Functions}

\begin{thm}~\label{mainbij} Fix a root set $R$ and let $\mathcal{MP} = \mathcal{MP_{R,G}}$ and $\mathcal{DC} = \mathcal{DC_{R,G}}$.  Define $\Omega : 
\mathcal{MP} \rightarrow \mathcal{DC}$ by $\Omega(f) = \Omega_f$ where

\[
\Omega_f(i) =
\left\{ \begin{array}{ll}
-\infty & \mbox{if $i \in R$.} \\
deg(i)-1-f(i) & \mbox{if $i \notin R$.}
\end{array} \right.
\]

\noindent Then $\Omega$ is a bijection, whose inverse $\Omega^{-1}: \mathcal{DC} \rightarrow \mathcal{MP}$, is given by $\Omega^{-1}(\mu) = 
\Omega^{-1}_{\mu}$ where

\[
\Omega^{-1}_{\mu}(i) =
\left\{ \begin{array}{ll}
\infty & \mbox{if $i \in R$.} \\
deg(i)-1-\mu(i) & \mbox{if $i \notin R$.}
\end{array} \right.
\]

\end{thm}

\begin{proof}

Let $f$ be any $G$-multiparking function.  First we show that $\Omega_f$ is a Dirichlet configuration.  This is trivial if $R=V(G)$, so assume $R 
\subset V(G)$.  As $\Omega_f(i) < deg(i)$ for every vertex $i$, $\Omega_f$ is stable.  By proposition \ref{recurrentcharacterization}, finding an 
avalanche connecting $(\Omega_f + \chi)$ to $\Omega_f$ is enough to show that $\Omega_f$ is recurrent.

Note that $(\Omega_f + \chi)(i) = deg(i) - 1 - f(i) + \chi(i)$ for every nonroot vertex $i$.  Therefore, a vertex $i$ in the configuration $(\Omega_f 
+ \chi)$ is ready if and only if $\chi(i) > f(i)$.  Since $f$ is a $G$-multiparking function, the set of all non-root vertices must have a 
well-behaved vertex, say $j$, and this implies $\chi(j) > f(j)$.  Hence, $(\Omega_f + \chi)$ is not stable.

Let $i_1, i_2, \ldots, i_n$ be a burning sequence for $f$ (in the sense of Lemma \ref{burning}), with $i_1, i_2, \ldots, i_k$ as the roots of $f$.  
We have just shown that there is a vertex that can be labelled $i_{k+1}$.  It is enough to show that if the vertices $i_{k+1}, \ldots, i_{l-1}$ can 
be fired, then $i_l$ can be fired.  Notice that for any $U \subseteq V(G)$, $deg(i_l) = \mathbb{O}_U(i_l) + \mathbb{I}_U(i_l)$.  So if $U = \{ i_l, 
i_{l+1}, \ldots, i_n \}$, then firing $i_{k+1}$ through $i_{l-1}$ sends exactly $\mathbb{O}_U(i_l) - \chi(i_l)$ chips to $i_l$.  So, $i_l$ will have 
at least $deg(i_l) - 1 - f(i_l) + \mathbb{O}_U(i_l)$ chips.  Since $f$ is a $G$-multiparking function, $f(i_l) < \mathbb{O}_U(i_l)$, so $deg(i_l) - 1 
- f(i_l) + \mathbb{O}_U(i_l) = deg(i_l) - 1 - (f(i_l) - \mathbb{O}_U(i_l)) \geq deg(i_l)$, and thus $i_l$ will be ready.  Hence, every non-root 
vertex in an avalanche beginning with $(\Omega_f + \chi)$ must be fired, and the throwing-out sequence specified is also a sequence in which the 
vertices can be fired.  (Note that although there may be several throwing-out sequences for $f$, they all yield the same final configuration.) Note 
that $\mu$ is a Dirichlet configuration if and only if a firing sequence exists, and this argument can be reversed to obtain a burning sequence, 
proving that this correspondence is surjective.

Finally, we must show that this sequence of firings beginning at $(\Omega_f + \chi)$ ends at $\Omega_f$.  If $i$ is any vertex, it loses $deg(i)$ 
chips when fired.  As its neighbors are fired, $i$ recovers exactly $deg(i) - \chi(i)$ chips, since the roots are not fired.  Thus, at the end of 
this avalanche, $i$ has exactly $deg(i) - 1 - f(i) + \chi(i) - deg(i) + (deg(i) - \chi(i)) = deg(i) -1 - f(i)$ chips, meaning that we end on the 
configuration $\Omega_f$.

Finally, it is obvious that $\Omega^{-1}$ is the inverse of $\Omega$.  \end{proof}

This result strengthens earlier work by Biggs (see Lemma 3(ii) in \cite{Bigg2}).  This simple bijection also provides information on the natural 
poset orders on the sets of $G$-multiparking functions and Dirichlet configurations with a given root set.  If $f$ is a $G$-multiparking function, it 
is immediate from the definition that any vertex function which is less than or equal to $f$ on each vertex is also a $G$-multiparking function.  
This determines a simple poset order on the $G$-multiparking functions.  Analogously, if $\mu$ is a Dirichlet configuration then any other 
configuration which is stable and greater than or equal to $\mu$ on every vertex is also Dirichlet.  Hence there is also a simple poset order on the 
Dirichlet configurations and the Hasse diagrams of these two posets are identical, except that one is upside-down.

\begin{cor} If $f$ and $g$ are $G$-multiparking functions, then $f \leq g$ (in the $G$-multiparking function poset order described above) if and only 
if $\Omega_f \geq \Omega_g$ (in the Dirichlet configuration poset order described above).

\end{cor}

Theorem \ref{mainbij} also suggests a burning-type algorithm for verifying that a configuration is Dirichlet for a given graph.

\begin{cor}~\label{critburning} A configuration $\mu$ on $G$ is Dirichlet $\Longleftrightarrow$ there exists a permutation $\pi \in S_n$ such that 
for every vertex $i$, either $\pi(i)$ is a root or $deg(\pi(i)) > \mu(\pi(i)) \geq \mathbb{I}_{U_i}(\pi(i))$, where $U_i := V(G) - \{ \pi(1), \ldots, 
\pi(i-1) \}$.

\end{cor}

\begin{proof}

By theorem \ref{mainbij}, $\mu$ is critical if and only if $f := \Omega^{-1}(\mu)$ is a $G$-multiparking function, and this is true if and only if 
there is a permutation $\pi \in S_n$ such that $0 \leq f(\pi(i)) < \mathbb{O}_{U_i}(\pi(i))$ for every nonroot vertex $i$.  But this is true if and 
only if
\begin{eqnarray*}
&& 0 \leq deg(\pi(i)) - 1 - \mu(\pi(i)) < \mathbb{O}_{U_i}(\pi(i)) \\
&\Longleftrightarrow& deg(\pi(i)) - 1 \geq \mu(\pi(i)) > deg(\pi(i)) - 1 - \mathbb{O}_{U_i}(\pi(i)) \\
&\Longleftrightarrow& deg(\pi(i)) > \mu(\pi(i)) \geq \mathbb{I}_{U_i}(\pi(i))
\end{eqnarray*}  \end{proof}

We will hereafter refer to the permutations in proposition \ref{critburning} as \emph{Dirichlet certificates for $\mu$}.  This proposition also helps 
us identify the avalanches connecting a Dirichlet configuration to itself.

\begin{prop}~\label{avalanchestructure} Let $\mu$ be a Dirichlet configuration and let $\pi \in S_n$.  $\pi$ is a Dirichlet certificate for $\mu$ 
$\Longleftrightarrow$ the avalanche determined by the firing sequence $\pi(1), \pi(2), \ldots, \pi(n)$ connects $\mu$ to itself.

\end{prop}

\begin{proof}

($\Longleftarrow$) Let $\alpha = (\mu=\mu_1, \mu_2, \ldots, \mu_n, \mu_1)$ be the avalanche determined by $\pi$ and suppose $\pi(i)$ is a nonroot.  
We must show that $deg(\pi(i)) > \mu_1(\pi(i)) \geq \mathbb{I}_{U_i}(\pi(i))$ for every such $i \leq n$.  By assumption, $\mu$ is a Dirichlet 
configuration, so $\mu$ is stable, and thus $deg(\pi(i)) > \mu_i(\pi(i))$ for every $i$.

The structure of $\alpha$ is that the vertices $\pi(1), \ldots, \pi(i-1)$ are fired, and after these firings we arrive at $\mu_i$.  Then 
$\mu_i(\pi(i)) - \mu_1(\pi(i)) = \mathbb{O}_{U_i}(\pi(i))$, since $U_i = V(G) - \{ \pi(1), \ldots, \pi(i-1) \}$.  Also, $\pi(i)$ is ready in $\mu_i$ 
and therefore $\mu_i(\pi(i)) \geq deg(\pi(i))$.  Thus 
\begin{eqnarray*}
\mathbb{O}_{U_i}(\pi(i)) &=& \mu_i(\pi(i)) - \mu_1(\pi(i)) \\
&\geq& deg(\pi(i)) - \mu_1(\pi(i)) \\
&=& \mathbb{O}_{U_i}(\pi(i)) + \mathbb{I}_{U_i}(\pi(i)) - \mu_1(\pi(i))
\end{eqnarray*}
Thus we have $\mu_1(\pi(i)) \geq \mathbb{I}_{U_i}(\pi(i))$, proving that $\pi$ is a Dirichlet certificate for $\mu_1$.

($\Longrightarrow$) If $\pi$ is a Dirichlet certificate, then $deg(\pi(i)) > \mu_1(\pi(i)) \geq \mathbb{I}_{U_i}(\pi(i))$ for every nonroot $\pi(i)$.  
Since $\mu = \mu_1$ is Dirichlet, it is stable, and thus only the roots can be fired.  Suppose $\pi(1), \ldots \pi(i-1)$ have been fired in that 
order.  Assuming $\pi(i)$ is not a root, $\mu_i(\pi(i)) = \mu_1(\pi(i)) + \mathbb{O}_{U_i}(\pi(i)) \geq \mathbb{I}_{U_i}(\pi(i)) + 
\mathbb{O}_{U_i}(\pi(i)) = deg(\pi(i))$, and so $\pi(i)$ is ready in $\mu_i$.  Thus, $\pi(1), \ldots, \pi(n)$ defines an avalanche.

It is clear that this avalanche connects $\mu_1$ to itself, since we begin at that configuration and every vertex is fired exactly once, meaning that 
the net change in chips at each vertex $i$ is $\sum_{j \neq i} e(i,j) - deg(i) = 0$.  \end{proof}

\end{section}


\begin{section}{Descending Traversals}

Let $G$ be as above, but connected and with a total ordering $<_E$ on the edge set $E(G)$ and $V(G)=[n]$.  Let $m = n + \# E(G)$.  Let $\Sigma = 
\Sigma(G) = (\sigma_1, \sigma_2, \ldots, \sigma_m)$ be a sequence of the edges and vertices of G in which each edge and vertex appears exactly once.  
Let $\Sigma^{\leq i} := (\sigma_1, \sigma_2, \ldots, \sigma_i)$ and $\Sigma^{\geq i} := (\sigma_i, \sigma_{i+1}, \ldots, \sigma_{m})$.  (Similarly, 
$\Sigma^{< i} := (\sigma_1, \sigma_2, \ldots, \sigma_{i-1})$ and $\Sigma^{> i} := (\sigma_{i+1}, \sigma_{i+2}, \ldots, \sigma_m)$.)  We define 
$\Sigma$ to be a {\em descending traversal on $G$} if it satisfies three conditions:

\begin{enumerate}

\item $\sigma_1$ is a vertex,

\item $\sigma_i$ ($i \neq 1$) a vertex $\Rightarrow$ $\sigma_{i-1}$ is an edge adjacent to $\sigma_i$,

\item $\sigma_i$ an edge $\Rightarrow$ it is adjacent to a vertex $\sigma_k$ with $k < i$ and $\sigma_i$ is maximal with respect to $<_E$ among all 
edges in $\Sigma^{\geq i}$ that are adjacent to some vertex in $\Sigma^{<i}$.

\end{enumerate}

This definition is due to Cori and LeBorgne \cite{CoLe}.  They provided explicit bijections from the descending traversals to the spanning trees and 
from the descending traversals to the critical configurations, and hence a bijection between these other two objects.

Now assume $G$ is the same as above, except not necessarily connected, and $R \subseteq V(G)$.  Let $\Sigma^*$ be a list of some vertices and edges 
of $G$ ($\Sigma^*$ contains no repetitions).  Let $\mathbb{E}(\Sigma^*)$ be the set of edges not in $\Sigma^*$ which are adjacent to a vertex in 
$\Sigma^*$.  We let $\Pi_i$ be the set of ordered pairs $(\Sigma^*,W)$ where $W \subseteq \mathbb{E}(\Sigma^*)$ and where not both of $W = \emptyset$ 
and $R \subseteq \Sigma^*$ is true.  A \emph{choice function} on $G$ is any function $\zeta$ from $\Pi_i$ to $E(G) \cup R$ such that

\[ \zeta(\Sigma^*,W) 
\left\{ \begin{array}{ll} 
\in W & \text{if } W \neq \emptyset \\ 
\in R - \Sigma^* & \text{if } W = \emptyset \text{ and } R \nsubseteq \Sigma^*
\end{array} \right. \]

Fix a choice function $\zeta$ and let $\Sigma = (\sigma_i)_{i=1}^{m}$ be a sequence containing each edge and vertex of $G$ exactly once.  We call 
$\Sigma$ a \emph{descending $R$-traversal on $G$}, where $R = \{ \sigma_{s_1}, \sigma_{s_2}, \ldots, \sigma_{s_k} \}$, such that each subsequence 
$S_i = (\sigma_{s_i}, \sigma_{s_i +1}, \ldots, \sigma_{s_{i+1}-1})$ of $\Sigma$ satisfies:

\begin{enumerate}

\item $\sigma_{s_i} = \zeta(\Sigma^{< s_i},\emptyset)$, where $\sigma_{s_i}$ is a root,

\item $\sigma_j \in S_i$, $j > s_i$, a vertex $\Rightarrow \sigma_{j-1}$ is an edge adjacent to $\sigma_j$,

\item $\sigma_j \in S_i$ an edge $\Rightarrow \sigma_j$ is adjacent to a vertex $\sigma_k$ with $k < j$ and $\sigma_j = \zeta(\Sigma^{\leq j-1}, 
\mathbb{E}(\Sigma^{\leq j-1}))$.

\end{enumerate}

Let $\mathcal{DT} = \mathcal{DT}_{R,G, \zeta}$ denote the set of descending $R$-traversals on $G$.  Note that the first condition and the requirement 
that $\Sigma$ can be partitioned into subsequences $S_i$ is not very restrictive.  To check that a subsequence is a descending $R$-traversal it is 
generally only necessary to confirm that the last two conditions hold.

If one defines $R := \{v_1\}$ and $\zeta$ to be the function that picks the largest-index edge available, then the descending $R$-traversals of $G$ 
are, in fact, just the descending traversals of $G$.

\begin{ex}~\label{someexs} We illustrate some descending $R$-traversals of $\Gamma$, for different $R$.  In all these examples, let 
$\zeta(\Sigma^*,W)$ be the largest-index edge in $W$ if $W \neq \emptyset$ and the smallest vertex in $R$ otherwise.

\begin{figure}[h]
\begin{picture}(0,0)
\put(200,-80){\includegraphics[width=4cm]{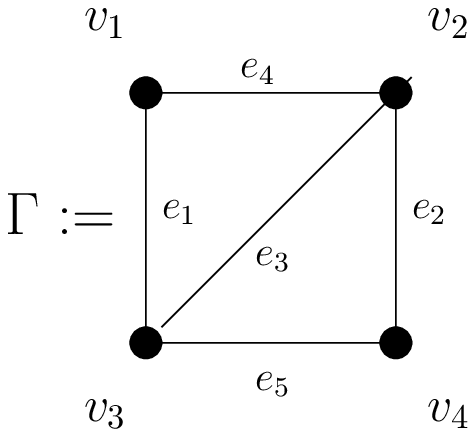}}
\end{picture}
\vskip 0.9 in
\end{figure}

\begin{enumerate}

\item Let $R = \{ v_1 \}$.  Then, $(v_1, e_4, v_2, e_3, e_2, v_4, e_5, v_3, e_1)$ and $(v_1, e_4, e_1, v_3, e_5, e_3, v_2, e_2, v_4)$ are descending 
$R$-traversals of $\Gamma$.

\item Let $R = \{ v_2, v_3 \}$.  Then, $(v_2, e_4, e_3, e_2, v_3, e_5, e_1, v_1, v_4)$ and $(v_2, e_4, e_3, e_2, v_4, e_5, v_3, e_1, v_1)$ are 
descending $R$-traversals of $\Gamma$.

\item Let $R = \{ v_1, v_2, v_4 \}$.  Then, $(v_1, e_4, e_1, v_2, e_3, v_3, e_5, e_2, v_4)$ is a descending $R$-traversal of $\Gamma$.  

\end{enumerate}

\end{ex}

Now suppose $\Sigma$ is a descending $R$-traversal, $R = \{ \sigma_{s_1}, \sigma_{s_2}, \ldots, \sigma_{s_k} \}$, and $\zeta$ is the choice function.  
With this as input, we define a function $f_{\Sigma}$ on $V(G)$ in the following way: \\

\noindent {\bf Algorithm A}

\begin{enumerate}

\item If $v = \sigma_{s_i}$ for some $i$, then set $f_{\Sigma}(v) = \infty$.

\item Otherwise, set $f_{\Sigma}(v)=j-1$, where $j$ is the number of edges adjacent to $v$ that precede $v$ in $\Sigma$.

\end{enumerate}

Note that if $v \notin R$, then by condition (2) of the definition of a descending $R$-traversal, it is preceded by an edge adjacent to it.  Thus, 
$f_{\Sigma}(v) = j-1 \geq 0$ and so $f_{\Sigma} : V(G) \rightarrow \mathbb{N} \cup \{ \infty \}$.

\begin{prop} $f_{\Sigma} \in \mathcal{MP}$ for any $\Sigma \in \mathcal{DT}$.

\end{prop}

\begin{proof} Let $f = f_{\Sigma}$ and let $\sigma_{v_1}, \sigma_{v_2}, \ldots, \sigma_{v_n}$ be the vertex subsequence of $\Sigma$.  We will show 
that this is a burning sequence for $f$, proving by Lemma \ref{burning} that $f$ is a $G$-multiparking function. (It is clear that $f$ has $k$ 
roots.)

First, note that $f(\sigma_{v_1}) = \infty$.  Let $U_i$ be the set of vertices in $\Sigma^{\leq i}$.  Now suppose $\sigma_{v_1}, \sigma_{v_2}, 
\ldots, \sigma_{v_{i-1}}$ are all either roots or well-behaved in $U_1, U_2, \ldots U_{i-1}$, respectively.  Suppose $\sigma_{v_i}$ is not a root.  
If $f(\sigma_{v_i}) = j-1$, then there are exactly $j$ edges adjacent to $\sigma_{v_i}$ and preceding it in $\Sigma$.  Each of these edges is 
preceded in $\Sigma$ by a vertex adjacent to it (note part $3$ of the definition of a descending $R$-traversal).  These vertices are among $\{ 
\sigma_{v_1}, \sigma_{v_2}, \ldots, \sigma_{v_{i-1}} \} = U_i$, and thus $0 \leq f(\sigma_{v_i}) = j-1 < j = \mathbb{O}_{U_i}(\sigma_{v_i})$.  Lemma 
\ref{burning} implies that $f \in \mathcal{MP}$.  \end{proof}

\begin{ex}~\label{psiex} Let $\zeta(\Sigma^*, W)$ be the largest-index edge if $W \neq \emptyset$ and the lowest-index vertex in $R - \Sigma^*$ 
otherwise.  Let $R = \{ v_1, v_4 \}$.  Below is a table listing some descending $R$-traversals of $\Gamma$ on the left-hand side and the 
corresponding (under algorithm A) $\Gamma$-multiparking functions on the right-hand side.  (The list of descending $R$-traversals is not exhaustive.)

\[
\begin{array}{lll}
(v_1, e_4, e_1, v_4, e_5, v_3, e_3, v_2, e_2) & \rightarrow & (\infty, 1, 1, \infty) \\
(v_1, e_4, e_1, v_4, e_5, v_3, e_3, e_2, v_2) & \rightarrow & (\infty, 2, 1, \infty) \\
(v_1, e_4, e_1, v_4, e_5, e_2, v_2, e_3, v_3) & \rightarrow & (\infty, 1, 2, \infty) \\
(v_1, e_4, v_2, e_3, e_2, v_4, e_5, e_1, v_3) & \multirow{2}{*}{\huge \}} & \multirow{2}{*}{$(\infty, 0, 2, \infty)$} \\
(v_1, e_4, v_2, e_3, e_2, e_1, v_4, e_5, v_3) & & \\
(v_1, e_4, e_1, v_3, e_5, v_4, e_3, v_2, e_2) & \multirow{2}{*}{\huge \}} & \multirow{2}{*}{$(\infty, 1, 0, \infty)$} \\
(v_1, e_4, e_1, v_3, e_5, e_3, v_2, e_2, v_4) & & \\
(v_1, e_4, e_1, v_3, e_5, v_4, e_3, e_2, v_2) & \multirow{2}{*}{\huge \}} & \multirow{2}{*}{$(\infty, 2, 0, \infty)$} \\
(v_1, e_4, e_1, v_3, e_5, e_3, v_4, e_2, v_2) & & \\
(v_1, e_4, v_2, e_3, e_2, v_4, e_5, v_3, e_1) & \multirow{2}{*}{\huge \}} & \multirow{2}{*}{$(\infty, 0, 1, \infty)$} \\
(v_1, e_4, v_2, e_3, e_2, e_1, v_3, e_5, v_4) & & \\
(v_1, e_4, v_2, e_3, v_3, e_5, e_2, e_1, v_4) & \multirow{2}{*}{\huge \}} & \multirow{2}{*}{$(\infty, 0, 0, \infty)$} \\
(v_1, e_4, v_2, e_3, v_3, e_5, e_2, v_4, e_1) & & 
\end{array} 
\]

\end{ex}

Lemma 1 of \cite{CoLe} states that if $(\sigma_i)_{i=1}^{m}$ and $(\tau_i)_{i=1}^{m}$ are descending traversals and $k$ is the minimal index at which 
they differ, then one of $\sigma_k$ and $\tau_k$ is an edge and the other is a vertex.  The example above shows that this is not necessarily true for 
descending $R$-traversals; $(v_1, e_4, e_1, v_4, e_5, e_2, v_2, e_3, v_3)$ and $(v_1, e_4, e_1, v_3, e_5, v_4, e_3, e_2, v_2)$ do not observe this 
property.

Let $\Psi = \Psi_{R,G,\zeta} : \mathcal{DT} \rightarrow \mathcal{MP}$ be defined by $\Psi(\Sigma) = f_{\Sigma}$.  The above example also illustrates 
that $\Psi$, as defined, is not generally injective.  We will now define, for each graph $G$, root set $R$, and choice function $\zeta$, a partition 
of $\mathcal{DT}$ over which $\Psi$ will turn out to be injective.

Let $f$ be any function from $V(G)$ to $\mathbb{N} \cup \{ \infty \}$ such that $f(i) = \infty$ if and only if $i \in R$.  We can consider 
$\Psi^{-1}(f)$, the (possibly empty) set of all descending $R$-traversals that are mapped to $f$.  It is then clear that $\mathcal{R} = 
\mathcal{R}_{R, G, \zeta} := \{ \Psi^{-1}(f) \; | \; \Psi^{-1}(f) \neq \emptyset \}$ is a partition of the set of descending $R$-traversals, where 
$\Psi^{-1}(f) = \{ \Sigma \in \mathcal{DT} \; | \; \Psi(\Sigma) = f \}$.  It is also clear that $\Psi$ is constant over each $\Psi^{-1}(f)$ in 
$\mathcal{R}$, and that $\Psi$ is injective when viewed as a function with $\mathcal{R}$ as its domain.  Throughout the rest of this paper, we will 
view $\Psi$ as a function from $\mathcal{R}$ to $\mathcal{MP}$.

Now we define an algorithm that will convert a a multiparking function to a descending $R$-traversal.\\

\noindent {\bf Algorithm B}

\begin{itemize}

\item {\bf Step 1: initial condition.} If $i=1$ then $\Sigma^{\leq i} := (\zeta(\emptyset,\emptyset))$.

\item {\bf Step 2: insert the next entry.} Suppose $i>1$.  If there exists a vertex $v \notin \Sigma^{\leq i-1}$ such that $\Sigma^{\leq i-1}$ 
contains exactly $f(v)+1$ edges adjacent to $v$, then $\Sigma^{\leq i} := < \Sigma^{\leq i-1},v >$.  If no such vertex exists, then $\Sigma^{\leq i} 
:= < \Sigma^{\leq i-1}, \zeta(\Sigma^{\leq i-1},\mathbb{E}(\Sigma^{\leq i-1}))> $.  Repeat this step until $i=m$.

\end{itemize}

\begin{ex} Recall the conditions in Example \ref{psiex}.  Below is a table listing all the $\Gamma$-multiparking functions on the right-hand side and 
the corresponding (under algorithm B) descending $R$-traversals of $\Gamma$ on the right-hand side.

\[
\begin{array}{lll}
(\infty, 1, 1, \infty) & \rightarrow & (v_1, e_4, e_1, v_4, e_5, v_3, e_3, v_2, e_2) \\
(\infty, 2, 1, \infty) & \rightarrow & (v_1, e_4, e_1, v_4, e_5, v_3, e_3, e_2, v_2) \\
(\infty, 1, 2, \infty) & \rightarrow & (v_1, e_4, e_1, v_4, e_5, e_2, v_2, e_3, v_3) \\
(\infty, 0, 2, \infty) & \rightarrow & (v_1, e_4, v_2, e_3, e_2, e_1, v_4, e_5, v_3) \\
(\infty, 1, 0, \infty) & \rightarrow & (v_1, e_4, e_1, v_3, e_5, e_3, v_2, e_2, v_4) \\
(\infty, 2, 0, \infty) & \rightarrow & (v_1, e_4, e_1, v_3, e_5, e_3, v_4, e_2, v_2) \\
(\infty, 0, 1, \infty) & \rightarrow & (v_1, e_4, v_2, e_3, e_2, e_1, v_3, e_5, v_4) \\
(\infty, 0, 0, \infty) & \rightarrow & (v_1, e_4, v_2, e_3, v_3, e_5, e_2, e_1, v_4) 
\end{array}
\]

%

\end{ex}

\begin{prop}~\label{inDT} $\Sigma^{\leq m} \in \mathcal{DT}$ for any $f \in \mathcal{MP}$.

\end{prop}

\begin{proof} 

We must first show that algorithm B can, in fact, always reach $\Sigma^{\leq m}$ if it acts on some $f \in \mathcal{MP}$.  Clearly $\Sigma^{\leq 
1}$ can be reached, so suppose $\Sigma^{\leq i} = (\sigma_j)_{j=1}^{i}$ can be reached for some $1 \leq i < m$.  There are two cases in which 
Algorithm B might fail to reach $\Sigma^{\leq m}$.

First, suppose that there are two vertices $v$ and $w$, neither in $\Sigma^{\leq i}$, such that there are exactly $f(v)+1$ and $f(w)+1$ edges in 
$\Sigma^{\leq i}$ that are adjacent to them respectively. (We may also assume, without loss of generality, that $i$ is the minimum index at which 
there is more than one vertex ready to be appended to $\Sigma^{\leq i}$.)  Note that the edge $\{ v,w \}$, if it exists, is not in $\Sigma^{\leq i}$; 
no such edge could be in $\mathbb{E}(\Sigma^{\leq j})$ for any $j \leq i$ since neither $v$ nor $w$ is in $\Sigma^{\leq j}$.  Therefore, in the 
sequence $\Sigma^{\leq i-1}, \Sigma^{\leq i-2}, \ldots, \Sigma^{\leq 1}$ there must be a $\Sigma^{\leq j}$ which contains exactly $f(v)+1$ edges 
adjacent to $v$ but fewer than $f(w)+1$ edges adjacent to $w$.  Thus, $v$ should have been added earlier and $i$ does not exist.

We must also show that there is no index $i$ for which $R \subseteq \Sigma^{\leq i}$ and $\mathbb{E}(\Sigma^{\leq i}) = \emptyset$.  Let $i$ be an 
index for which there is no vertex $v \notin \Sigma^{\leq i}$ adjacent to exactly $f(v)+1$ edges in $\Sigma^{\leq i}$.  Assume $R \subseteq 
\Sigma^{\leq i}$.  Clearly, if $\Sigma^{\leq i}$ contains $V(G)$, then $\mathbb{E}(\Sigma^{\leq i})$ cannot be empty unless $i=m$.  So let $U$ be the 
set of vertices not in $\Sigma^{\leq i}$.  Since $f \in \mathcal{MP}$ and there is no root in $U$, this set must have a well-behaved vertex.  That 
is, there is a vertex $v \in U$ such that $0 \leq f(v) < \mathbb{O}_U(v)$.  In particular, $0 < \mathbb{O}_U(v) \leq \# \mathbb{E}(\Sigma^{\leq i})$.

It is clear, from the construction of Algorithm B, that $\Sigma^{\leq m}$ satisfies the last two conditions in the definition of a descending 
$R$-traversal.  \end{proof}

Let $\Phi = \Phi_{R, G, \zeta} : \mathcal{MP} \rightarrow \mathcal{DT}$ be defined by $\Phi(f) = \Sigma^{\leq m}$.

\begin{prop}~\label{phiinj} $\Phi$ is injective.

\end{prop}

\begin{proof} Let $f$ and $g$ be different functions in $\mathcal{MP}$, and $\Phi(f) = \Sigma_{f}$ and $\Phi(g) = \Sigma_{g}$.  Since $f$ and $g$ are 
different, there is a vertex $v$ at which $f(v) < g(v)$ ($v$ is not a root, since $f$ and $g$ have the same root set).  There is an index $i$ at 
which $v$ appears in $\Sigma_{f}$.  This means $\Sigma_{f}^{\leq i} = < \Sigma_{f}^{\leq i-1}, v >$, but then either $\Sigma_f^{\leq i-1} \neq 
\Sigma_g^{\leq i-1}$ or $\Sigma_{g}^{\leq i} \neq < \Sigma_{g}^{\leq i-1}, v >$ and so $\Phi(f) \neq \Phi(g)$.  \end{proof}
	
\begin{prop} $\Psi(\Phi(f)) = f$ for any $f \in \mathcal{MP}$.

\end{prop} 

\begin{proof} It is enough to show that $\Phi(f) \in \Psi^{-1}(f)$ for any $f \in \mathcal{MP}$.  Suppose $\Phi(f) = (\sigma_i)_{i=1}^{m}$.  Note 
that $\sigma_i \in R$ if and only if $\sigma_i = \zeta(\Sigma^{\leq i-1}, \emptyset)$ (where $\Sigma^{\leq i-1} = \emptyset$ if $i=1$) and this is 
true if and only if $f(\sigma_i) = \infty$.  Therefore $\Psi(\Phi(f))|_{\sigma_i} = \infty = f(\sigma_i)$.  If $f(\sigma_i) = a$ for some vertex 
$\sigma_i$ then $\Sigma^{\leq i-1}$ contains exactly $a + 1$ edges adjacent to $\sigma_i$.  Therefore, $\Psi(\Phi(f))|_{\sigma_i} = a = f(\sigma_i)$.  
So, $\Psi$ maps $\Phi(f)$ to $f$ and thus $\Phi(f) \in \Psi^{-1}(f)$.  \end{proof}

\end{section}


\begin{section}*{Acknowledgements}

I extend my sincere gratitude to Professors Rob Ellis and Catherine Yan for their many helpful comments.

\end{section}



\end{document}